# MULTIVARIATE SPACINGS BASED ON DATA DEPTH: I. CONSTRUCTION OF NONPARAMETRIC MULTIVARIATE TOLERANCE REGIONS[1]

By Jun Li and Regina Y. Liu

*University of California at Riverside and Rutgers University*

This paper introduces and studies multivariate spacings. The spacings are developed using the order statistics derived from data depth. Specifically, the spacing between two consecutive order statistics is the region which bridges the two order statistics, in the sense that the region contains all the points whose depth values fall between the depth values of the two consecutive order statistics. These multivariate spacings can be viewed as a data-driven realization of the so-called "statistically equivalent blocks." These spacings assume a form of center-outward layers of "shells" ("rings" in the two-dimensional case), where the shapes of the shells follow closely the underlying probabilistic geometry. The properties and applications of these spacings are studied. In particular, the spacings are used to construct tolerance regions. The construction of tolerance regions is nonparametric and completely data driven, and the resulting tolerance region reflects the true geometry of the underlying distribution. This is different from most existing approaches which require that the shape of the tolerance region be specified in advance. The proposed tolerance regions are shown to meet the prescribed specifications, in terms of $\beta$-*content* and $\beta$-*expectation*. They are also asymptotically minimal under elliptical distributions. Finally, a simulation and comparison study on the proposed tolerance regions is presented.

**1. Introduction.** The term "spacings" in statistics generally refers to either the intervals (or gaps) between two consecutive order statistics or the lengths of these intervals. Spacings have been used extensively in probability and statistics, especially in the areas of distributional characterization,

Received January 2007; revised May 2007.

[1]Supported in part by grants from the National Science Foundation, the National Security Agency, and the Federal Aviation Administration.

extreme value theory and nonparametric inference. There is a rich literature on the theory and applications of spacings. The excellent treatise by Pyke in [22] as well as the references therein (e.g., [5] and [27]) and thereafter (e.g., [2, 4, 9, 12, 28]) all attest to the importance of spacings. In his paper [22] Pyke wrote,

> Perhaps the most significant restrictions of this paper has been our concern with one-dimensional spacings. There are many applications in which samples are drawn from two- or even three-dimensional space and for which it is important to study the spacings of the observations.

Although research on spacings has continued, his call for multivariate spacings has remained largely unanswered. The main difficulty in generalizing the univariate spacings to multivariate settings is the lack of suitable ordering schemes for multivariate observations. This paper has two goals. First, we introduce multivariate spacings using the multivariate ordering derived from the notion of data depth. Second, as an application, we apply the proposed multivariate spacings to construct nonparametric tolerance regions.

The paper is organized as follows. Section 2 is devoted to the development of multivariate spacings. We begin with a brief review of the univariate spacings and of some of their properties, as well as a brief description of the subject of data depth and the corresponding depth ordering of multivariate data. Note that the depth ordering is from the center-outward rather than the usual univariate linear ordering from the smallest to the largest. For any two consecutive depth order statistics, we define the spacing between them as the region that contains all the points in the sample space whose depth values fall between the depth values of the two order statistics. The multivariate spacings are the collection of these regions formed by all pairs of consecutive order statistics. These regions generally appear as center-outward layers of "shells" ("rings" in $\Re^2$), and the shapes of the shells follow closely the probabilistic geometry of the underlying distribution. In Section 3 we first provide a review of tolerance intervals for univariate data as well as the existing approaches for obtaining multivariate tolerance regions. We then describe the construction of nonparametric tolerance regions using the proposed multivariate spacings, and investigate the properties of the proposed tolerance regions. Specifically, we show that these tolerance regions: (1) meet the prescribed specifications in terms of $\beta$-*content* and $\beta$-*expectation*, and (2) are asymptotically minimal under a certain class of distributions which includes the elliptical family. The formation of our tolerance region is completely data driven and nonparametric, and the resulting tolerance region has the desirable property of reflecting accurately the underlying probabilistic geometry. In other words, the shape of our proposed tolerance regions is automatically determined by the given data, and does not need to be specified in advance. Most existing approaches require pre-specification of the shape, which can be considered arbitrary or subjective. It is also worth noting that our tolerance region is always connected, which is more suitable in applications such as quality control. Section 4 contains a simulation study and some comparisons with other existing tolerance regions. It confirms several desirable features of our approach. Section 5 contains some concluding remarks. Most technical proofs are collected in the Appendix.

We also observe in Section 3.1 that in using our multivariate spacings to construct tolerance regions, we have in effect argued that our multivariate spacings are an ideal realization of the so-called "statistically equivalent blocks." This is because the realization of our multivariate spacings and their shapes are entirely data driven. Statistically equivalent blocks had been considered by Tukey in [24] and several follow-up papers (see, e.g., [10]) as possible building blocks for the construction for tolerance regions or tools for characterizing distributions. However, these papers again all need to pre-specify the shapes (e.g. rectangles or circles in $\Re^2$) of the blocks.

**2. Multivariate spacings derived from data depth.** We begin with a brief review of the notion of data depth and its corresponding multivariate ordering. This multivariate ordering naturally leads to our multivariate spacings.

2.1. *Data depth and center-outward ordering of multivariate data.* A *data depth* is a measure of "depth" of a given point with respect to a multivariate data cloud or its underlying distribution, and it gives rise to a natural *center-outward ordering* of the points in a multivariate sample. Although the actual depth value has been used widely to develop robust multivariate inference, the depth-ordering is less understood and still underutilized. Existing notions of data depth include: Mahalanobis depth ([20]), half-space depth ([14, 25]), simplicial depth ([16]), projection depth ([7, 8, 23, 30]), etc. More discussion on different notions of data depth can be found in [17, 31].

To help facilitate the coming exposition of multivariate spacings, we use the simplicial depth to illustrate the general concept of data depth and its corresponding center-outward ordering. Let $\{X_1, \ldots, X_n\}$ be a random sample from the distribution $F(\cdot) \in \Re^p$, $p \geq 2$. Consider the bivariate setting, $p = 2$. Let $\triangle(a, b, c)$ denote the triangle with vertices $a$, $b$ and $c$. Let $I(\cdot)$ be the indicator function, that is $I(A) = 1$ (or 0) if $A$ occurs (or not). For the given sample $\{X_1, \ldots, X_n\}$, the sample simplicial depth of $x$ is defined as

$$(2.1) \qquad D_{F_n}(x) = \binom{n}{3}^{-1} \sum_{(*)} I(x \in \triangle(X_{i_1}, X_{i_2}, X_{i_3})),$$

which is the fraction of the triangles generated from the sample that contain the point $x$. Here $(*)$ runs over all possible triplets of $\{X_1, \ldots, X_n\}$. A larger value of $D_{F_n}(x)$ indicates that $x$ falls in more triangles generated from the sample, and thus lies deeper within the data cloud.

The above can be generalized to dimension $p$ by counting simplices rather than triangles, that is

$$(2.2) \quad D_{F_n}(x) = \binom{n}{p+1}^{-1} \sum_{(*)} I(x \in s[X_{i_1}, \ldots, X_{i_{p+1}}]),$$

where $(*)$ runs over all possible subsets of $\{X_1, \ldots, X_n\}$ of size $(p+1)$. Here $s[X_{i_1}, \ldots, X_{i_{p+1}}]$ is the closed simplex whose vertices are $\{X_{i_1}, \ldots, X_{i_{p+1}}\}$.

If $F$ is given, the simplicial depth of $x$ w.r.t. to $F$ is defined as $D_F(x) = P_F\{x \in s[X_1, \ldots, X_{p+1}]\}$, where $X_1, \ldots, X_{p+1}$ are $(p+1)$ random observations from $F$. $D_F(x)$ measures how "deep" $x$ is w.r.t. $F$, and $D_{F_n}(x)$ in (2.2) is a sample estimate of $D_F(x)$. A fuller motivation together with the key properties of data depth can be found in [16]. In particular, it is shown that $D_F(\cdot)$ is affine invariant, and that $D_{F_n}(\cdot)$ converges uniformly and strongly to $D_F(\cdot)$. The affine invariance ensures that our proposed spacings and inference methods are coordinate free, and the convergence of $D_{F_n}$ to $D_F$ allows us to approximate $D_F(\cdot)$ by $D_{F_n}(\cdot)$ if $F$ is unknown.

For the given sample $\{X_1, X_2, \ldots, X_n\}$, we calculate the depth values $D_{F_n}(X_i)$'s and then order the $X_i$'s according to their descending depth value. Denoting by $X_{[j]}$ the sample point associated with the $j$th largest depth value, we then obtain the sequence $\{X_{[1]}, X_{[2]}, \ldots, X_{[n]}\}$ which is the depth order statistics of $X_i$'s, with $X_{[1]}$ being the *deepest* point, and $X_{[n]}$ most outlying. Here, a larger order is associated with a more outlying position w.r.t. the underlying distribution. Note that the order statistics derived from depth are different from the usual order statistics in the univariate case, since the latter are ordered from the smallest sample point to the largest, while the former is from the *middle* sample point and moves outward in all directions. Figure 1 helps demonstrate this feature of the depth ordering. The two plots show two random samples, each of size 500, drawn respectively from the standard bivariate normal and bivariate exponential distributions. For each plot, the "+" marks the deepest point, and the most inner convex hull encloses the deepest 20% of the sample points. The convex hull expands outward to enclose the next deepest 20% by each expansion. Those convex hulls determined by the decreasing depth value are nested, a feature indicating that the depth ordering is from the center outward. Note that the shape of the depth contours in those plots clearly reflects the underlying probabilistic geometry, relatively spherical in the normal case and fanning upper-right triangularly in the exponential case. The nested shells-like depth contours in Figure 1 also help illustrate the features of the multivariate spacings in Section 2.3.

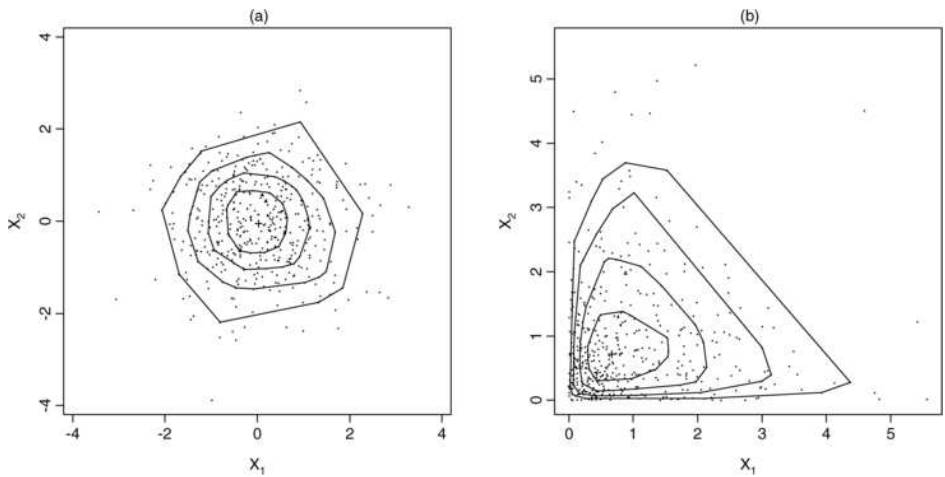

FIG. 1. *Depth contours for:* (a) *bivariate normal sample;* (b) *bivariate exponential sample.*

We give the definition of Mahalanobis depth here, since it is also used in the simulation study later in Section 4.

DEFINITION 2.1. The *Mahalanobis depth* ([20]) at $x$ with respect to $F$ is defined as

$$mD_F(x) = [1 + (x - \mu_F)\Sigma_F^{-1}(x - \mu_F)']^{-1},$$

where $\mu_F$ and $\Sigma_F$ are the mean vector and dispersion matrix of $F$, respectively. The sample version of the Mahalanobis depth is obtained by replacing $\mu_F$ and $\Sigma_F$ with their sample estimates.

Different notions of depth are capable of capturing different aspects of the probabilistic geometry, and may lead to different ordering schemes. However, all the depth orderings are essentially from the center outward. We note that all the depths aforementioned are affine invariant, and so are their resulting orderings. The affine invariance is a desirable feature for the construction of multivariate spacings later in Section 2.3.

Note that geometric depths such as the half-space and the simplicial depths are completely nonparametric and moment-free, and they capture well the underlying probabilistic geometry of the data. Although the Mahalanobis depth captures less well the underlying geometry unless the geometry happens to be elliptical, it is computationally more feasible than geometric depths. Under elliptical distributions, the two geometric depths capture fairly well the elliptical structure in the large sample case and are close competitors to the Mahalanobis depth. Between the two geometric depths, the simplicial depth provides a finer ordering and produces less ties than the half-space depth. This point has been observed in [18].

For convenience, we will use the notation $D(\cdot)$ to express any valid notion of depth, unless a particular notion is to be emphasized.

Before we use depth order statistics to formulate multivariate spacings, we review the univariate spacings and some of their properties.

2.2. *Univariate spacings.* Let $X_1, X_2, \ldots, X_n$ be a random sample from a univariate continuous distribution $F$ which has the support $(a, b)$. Denote by $\{X_{[1]}, X_{[2]}, \ldots, X_{[n]}\}$ the order statistics of $X_i$'s, namely $X_{[1]} \leq X_{[2]} \leq \cdots \leq X_{[n]}$. Note that we will avoid introducing additional messy notation for differentiating the univariate setting from the multivariate one by using the same notation $X_{[j]}$ throughout the paper to indicate the $j$th order statistic of the sample $X_i$'s. It will generally be clear from the context whether the notation is intended for the univariate or for the multivariate setting. If needed, the phrase "univariate ordering" or "depth ordering" will be used to emphasize the univariate or the multivariate ordering.

Given the order statistics $X_{[1]} \leq X_{[2]} \leq \cdots \leq X_{[n]}$, the univariate spacings of the sample refer to the intervals $L_i = (X_{[i-1]}, X_{[i]}]$, $i = 1, \ldots, n+1$, with $X_{[0]} = a$ and $X_{[n+1]} = b$, or their lengths $D_i = X_{[i]} - X_{[i-1]}$. For convenience, we proceed to discuss the spacings by assuming that $F$ follows the uniform distribution on $(0, 1)$ [denoted by $F \sim U(0,1)$], since the probability integral transformation $F(X)$ transforms the given sample into a sample from $U(0,1)$. If $F \sim U(0,1)$, then:

(i) $D_1 + D_2 + \cdots + D_{n+1} = 1$, and
(ii) the density function of $(D_1, D_2, \ldots, D_{n+1})$ is

$$f(d_1, d_2, \ldots, d_{n+1}) = \begin{cases} n!, & \text{if } d_i \geq 0 \text{ and } d_1 + d_2 + \cdots + d_{n+1} = 1, \\ 0, & \text{otherwise.} \end{cases}$$

Thus the density function $f$ is completely symmetrical in its arguments.

[21] and [22] have observed that the uniform spacings $(D_1, D_2, \ldots, D_{n+1})$ can be viewed as exponential random variables proportional to their sum. Specifically, assume that $\{U_1, U_2, \ldots, U_{n+1}\}$ is a random sample from the exponential distribution with mean 1 [denoted as $Exp(1)$], and let

$$S = U_1 + U_2 + \cdots + U_{n+1} \quad \text{and} \quad W_i = U_i/S, \qquad i = 1, \ldots, n+1.$$

Then, $(W_1, W_2, \ldots, W_{n+1})$ and $(D_1, D_2, \ldots, D_{n+1})$ are identically distributed.

A similar property will appear during our formulation of multivariate spacings later.

2.3. *Multivariate spacings.* The main difficulty in extending the univariate spacings to higher dimensions lies in the lack of proper ordering of the multivariate data. Applying the center-outward ordering induced from data depth, the multivariate spacings can be defined as follows.

Let $X_1, \ldots, X_n$ be a random sample from a continuous distribution $F$ in $\Re^p$, $p \geq 2$. For a given data depth $D(\cdot)$, we calculate all $D_F(X_i)$'s, and obtain the depth order statistics $X_{[1]}, \ldots, X_{[n]}$ in descending depth values. Let $Z_i = D_F(X_i)$, and $Z^{[i]} = D_F(X_{[i]})$, for $i = 1, \ldots, n$. Note that $Z^{[1]} \geq \cdots \geq Z^{[n]}$, which are reverse univariate order statistics of $Z_i$'s. The matching indices in $Z^{[1]} \geq \cdots \geq Z^{[n]}$ and $X_{[1]}, \ldots, X_{[n]}$ are useful for tracking depth order statistics with their depth values in defining multivariate spacings and tolerance regions later. We now define the multivariate spacings as follows,

(2.3) $\quad MS_i = \{X : Z^{[i-1]} \geq D_F(X) > Z^{[i]}\}, \qquad i = 1, \ldots, n+1,$

with $Z^{[0]} = \sup_x \{D_F(x)\}$ and $Z^{[n+1]} = 0$. The corresponding sample multivariate spacings are

(2.4)
$$\widehat{MS}_i = \{X : \hat{Z}^{[i-1]} \geq D_{F_n}(X) > \hat{Z}^{[i]}\}, \qquad i = 1, \ldots, n, \quad \text{and}$$
$$\widehat{MS}_{n+1} = \{X : D_{F_n}(X) \leq \hat{Z}^{[n]}\},$$

where $\hat{Z}^{[0]} = \sup_x \{D_{F_n}(x)\}$, and $\hat{Z}^{[1]} \geq \cdots \geq \hat{Z}^{[n]}$ are the reverse order statistics of $\hat{Z}_i = D_{F_n}(X_i)$, $i = 1, \ldots, n$.

Note that the multivariate spacings here define the "gap" between two consecutive depth order statistics as the shell-shape region bridging the two order statistics, generalizing the interval linking the two consecutive order statistics in the univariate spacings. Consequently, the multivariate spacings derived from depth order statistics are center-outward layers of "shells." Figure 2 illustrates an example of multivariate spacings determined by a random sample of size five drawn from the bivariate normal distribution with mean $(0,0)$ and covariance matrix $\begin{pmatrix} 3 & 1 \\ 1 & 1 \end{pmatrix}$. The five data points are denoted by circles in the plot. The Mahalanobis depth is used to calculate depth values. The multivariate spacings include six regions, five center-outward layered shells and the outmost region. Note that the shells clearly reflect the elliptical shape of the underlying distribution. Plots of the multivariate spacings for the standard bivariate normal and exponential samples using the simplicial depth show layered shells with shapes similar to those of Figure 1. Again, the shape of shells reflects the underlying geometric features.

Next, we observe a useful property regarding the coverage probabilities of the proposed multivariate spacings.

THEOREM 2.1. *Let $X_1, \ldots, X_n$ be a random sample from $F \in \Re^p$. Assume that the notion of data depth used in deriving the multivariate spacings (2.3) is affine invariant. Then, the coverage probabilities of these multivariate spacings, namely $\{P_F(MS_1), \ldots, P_F(MS_{n+1})\}$, follow the same distribution as the univariate spacings $\{D_1, \ldots, D_{n+1}\}$.*

PROOF. Let $Z_i = D_F(X_i)$ and $T_i = P_F(X : D_F(X) > Z_i)$, for $i = 1, \ldots, n$. Then $T_i$'s can be considered as a random sample drawn from $U[0,1]$, as seen in [19]. Let $T_{[1]} \leq \cdots \leq T_{[n]}$ be the order statistics of $T_i$'s. It is clear that $T_{[i]} = P_F(X : D_F(X) > Z^{[i]})$. Therefore $P_F(MS_i) = T_{[i]} - T_{[i-1]}$, where $T_{[0]} = 0$ and $T_{[n+1]} = 1$, and thus the theorem follows. $\square$

**3. Tolerance region based on multivariate spacings.** A confidence interval is used to provide an interval estimate for a parameter of interest with a stated confidence level. In production processes or quality control, it is customary to seek an interval that covers a certain proportion of the process distribution with a stated confidence as an assurance for meeting the required product specification. Intervals which fulfill this need are called tolerance intervals. In many practical situations, the quality of a product is specified by multiple characteristics of the product. To ensure the specifications of those multiple characteristics simultaneously, multivariate tolerance regions are needed. Tolerance intervals and regions are integral parts of applications in reliability theory and quality control. They allow the control of intended proportions of productions to meet the specified requirements. A high percentage of the production outside this interval (or region) will result in a high loss or rework rate. Before we describe our proposed construction of tolerance regions, we briefly review the literature of tolerance intervals and regions.

If the underlying process distribution is known, from either the design of experiment or the knowledge gained over long experience, the tolerance intervals or regions usually can be established. For example, if the sample is drawn from $N(\mu, \sigma)$, a normal distribution with the known mean $\mu$ and variance $\sigma^2$, and if we define tolerance intervals as those which contains $100\beta\%$ of the underlying distribution, then the shortest tolerance interval is simply $(\mu - z_{(1-\beta)/2}\sigma, \mu + z_{(1-\beta)/2}\sigma]$. Here $z_{(1-\beta)/2}$ is the upper $(1-\beta)/2$th quantile of the standard normal distribution. The constant $\beta$ is referred to as the *tolerance level*. This development can be extended in a straightforward manner to the setting of a $p$-dimensional normal distribution with mean vector $\mu$ and covariance matrix $\Sigma$ [denoted as $\mathbf{N}(\mu, \Sigma)$]. In this case, the corresponding smallest tolerance region can be constructed as an ellipsoid. It follows the elliptical level sets of the underlying multivariate normal distribution and satisfies

$$Q_{\mu,\Sigma}(t) \equiv \{X : (X-\mu)^T \Sigma^{-1}(X-\mu) \leq t\},$$

where $t$ is the solution of the equation

$$\int \cdots \int_{r^2 \leq t} \sqrt{2\pi} e^{-(1/2)r^2} r^{p-1} \prod_{i=2}^{p-1} \sin^{p-i} \theta_{i-1} \, dr \, d\theta_1 \cdots d\theta_{p-1} = \beta.$$

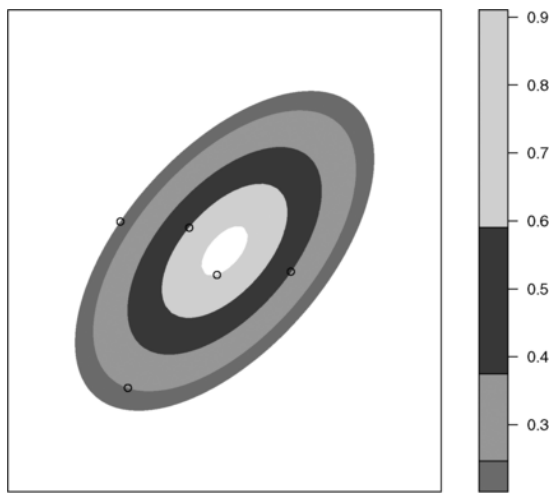

FIG. 2. *Multivariate spacings for a bivariate normal sample.*

If the distribution or its parameters are unknown, the following two definitions of tolerance regions have been considered and accepted as standard definitions, see [11] for example. Again, let $X_1, \ldots, X_n$ be a random sample from $F \in \Re^p$, $p \geq 1$.

DEFINITION 3.1. $T(X_1, \ldots, X_n)$ is called a $\beta$-*content* tolerance interval (or region) at confidence level $\gamma$ if

(3.1) $$P(P_F(T(X_1, \ldots, X_n)) \geq \beta) = \gamma.$$

DEFINITION 3.2. The region $T(X_1, \ldots, X_n)$ is called a $\beta$-*expectation* tolerance interval (or region) if

(3.2) $$E(P_F(T(X_1, \ldots, X_n))) = \beta.$$

In the univariate case, if the normality assumption holds but the parameters are unknown, a tolerance interval can be constructed by

(3.3) $$(\bar{X} - cS, \bar{X} + cS],$$

where $\bar{X}$ and $S$ are respectively the sample mean and standard. If Definition 3.1 is followed, [15] shows that $c$ can be approximated by $\sqrt{\frac{(n-1)(1+1/n)z_{(1-\beta)/2}^2}{\chi_{\gamma,n-1}^2}}$. Here $\chi_{\gamma,n-1}^2$ is the $(1-\gamma)$th quantile of the chi-square distribution with degree of freedom $(n-1)$. When the normality assumption is uncertain, Wilks (in [29]) proposed to use the order statistics, $X_{[1]} < \cdots < X_{[n]}$, to construct the following nonparametric tolerance interval:

$$(3.4) \qquad T(X_1, \ldots, X_n) = (X_{[r]}, X_{[n-r+1]}],$$

where $r$ is a positive integer and $r < (n+1)/2$. It has been shown that the coverage probability of this tolerance region, namely $P_F((X_{[r]}, X_{[n-r+1]}])$, follows a *Beta* distribution with parameters $(n-2r+1)$ and $2r$, denoted as $Beta(n-2r+1, 2r)$. Based on this observation, $r$ can be chosen to satisfy

$$(3.5) \qquad P(Beta(n-2r+1, 2r) \geq \beta) = \gamma$$

or

$$(3.6) \qquad E(Beta(n-2r+1, 2r)) = \beta$$

to meet the requirement in Definitions 3.1 or 3.2. Note that the tolerance interval in (3.4) is "symmetric" around the observed center point in the sense that the interval excludes an equal number of sample points from both tails. Wald in [26] considered a generalization of this symmetric tolerance interval, namely $(X_{[s]}, X_{[t]}]$, where $1 \leq s \leq t \leq n$. Clearly, this includes Wilks' interval as a special case, if $s = r$ and $t = n - r + 1$. Since the coverage probability of $(X_{[s]}, X_{[t]}]$ can be shown to follow $Beta(t-s, n-t+s+1)$, the desired tolerance interval for Definitions 3.1. or 3.2 can be obtained by choosing $s$ and $t$ as the solutions of

$$(3.7) \qquad P(Beta(t-s, n-t+s+1) \geq \beta) = \gamma$$

or

$$(3.8) \qquad E(Beta(t-s, n-t+s+1)) = \beta.$$

Note that the solution for (3.7) or (3.8) may not be unique. Different applications may impose different additional desirable properties and thus constraints on the choice of $s$ and $t$. One intuitively appealing and desirable property is that the tolerance interval (or region) be minimal.

To achieve the minimal nonparametric tolerance interval, Charterjee and Patra in [3] proposed a large-sample approach based on nonparametric density estimation, which yields asymptotically minimal tolerance intervals. The performance of this approach depends heavily on the methods used for density estimation and smoothing. Moreover, this approach tends to be overly conservative, as observed in [6]. When the underlying distribution is multimodal, the tolerance interval obtained by this approach may be the union of disjoint intervals, which is not desirable in practice.

In the multivariate case, when $F$ is unknown, there have been efforts to develop nonparametric multivariate tolerance regions. For example, Wald in [26] extended Wilks' approach for constructing the tolerance intervals in the univariate case to the multivariate case by sequentially adapting it for each coordinate. Under this method, the shape of the resulting tolerance region would be limited to the hyperrectangles (or rectangular blocks) with faces parallel to the coordinate hyperplanes. Tukey in [24] generalized Wald's approach to any desired shape by introducing the concept of "statistically equivalent blocks." However, the construction of the statistical equivalent blocks here requires choosing a priori an ordering function and thus can be somewhat arbitrary. Moreover, the shape of the constructed tolerance region based on this predetermined ordering function may be difficult to interpret or implement in practice. More discussion on statistically equivalent blocks is given later in Remark 3.1.

Chatterjee and Patra's approach for constructing asymptotically minimal tolerance intervals based on nonparametric density estimation is also applicable to the multivariate case, although it has the same drawbacks mentioned in the univariate case. Recently, using empirical process theory, Di Bucchianico, Einmahl and Mushkudiani [6] succeeded in developing an important new method for constructing the smallest nonparametric multivariate tolerance regions. Although this method possesses several desirable properties, it still has the following potential drawbacks: (i) it requires pre-specifying the shape of the tolerance region; (ii) the obtained tolerance region may not represent well the underlying geometry of the data; and (iii) the obtained region may not be connected. Finally, the computation involved in finding this smallest tolerance region can be quite intensive.

The above review of nonparametric multivariate tolerance regions shows that almost all existing approaches require specifying in advance the shape of the region. Most shapes specified, such as rectangular or elliptical, seem arbitrary and chosen mainly for mathematical convenience. If the shape is not chosen properly, these approaches may lead to gross misrepresentation of the underlying geometry of the data.

Recall that the nonparametric tolerance interval proposed in the Wilks approach [29] has the form $T(X_1, \ldots, X_n) = (X_{[r]}, X_{[n-r+1]}]$. From the point of view of spacings, the Wilks' tolerance interval can be easily seen as the union of some suitable number of the univariate spacings,

$$(3.9) \qquad T(X_1, \ldots, X_n) = (X_{[r]}, X_{[n-r+1]}] = \bigcup_{i=r+1}^{n-r+1} L_i.$$

Similarly, the proposed multivariate spacings derived in Section 2.3 can be used to form tolerance regions in multivariate settings. We now give details on such constructions, and discuss their properties.

3.1. *Properties of tolerance regions*: $F$ *is known.* Consider the case where $F \in \Re^p$ is known, $p \geq 2$. Recall that $X_{[1]}, \ldots, X_{[n]}$ denote the depth order statistics of the sample $X_i$'s and that $Z^{[1]}, \ldots, Z^{[n]}$ are their corresponding depth values. Recall also that $Z_i = D_F(X_i)$ and $Z^{[1]} \geq \cdots \geq Z^{[n]}$. Then we propose to form the tolerance region as the union of a suitable number of the inner spacings, which can be expressed as follows:

$$(3.10) \qquad O_{Z^{[r_n]}} = \bigcup_{i=1}^{r_n} MS_i \equiv \{X : D_F(X) > Z^{[r_n]}\},$$

for a suitably chosen $r_n$. Here $MS_i$, is the $i$th spacing, as defined in (2.3).

Applying Theorem 2.1, the distribution of the coverage probability of the above tolerance region can be determined immediately, as shown in the following theorem.

THEOREM 3.1. *The distribution of $P_F(O_{Z^{[r_n]}})$, the coverage probability of the tolerance region defined in* (3.10), *follows $Beta(r_n, n+1-r_n)$.*

PROOF. Clearly, $O_{Z^{[r_n]}} = \bigcup_{i=1}^{r_n} MS_i$. It follows from Theorem 2.1 that $\sum_{i=1}^{r_n} P_F(MS_i)$ and $\sum_{i=1}^{r_n} D_i$ are identically distributed. Here $D_i, i = 1, \ldots, n+1$, are the uniform spacings. Recalling the construction of the uniform spacings using exponential random variables given in Section 2.2, we see then

$$P_F(O_{Z^{[r_n]}}), \qquad \sum_{i=1}^{r_n} D_i, \quad \text{and} \quad \sum_{i=1}^{r_n} U_i \Big/ \sum_{j=1}^{n+1} U_j,$$

all have the same distribution. Here $U_1, U_2, \ldots, U_{n+1}$ are i.i.d. exponential random variables with mean 1, $Exp(1)$. Since $Exp(1)$ can also be viewed as the Gamma random variable $Gamma(1,1)$, $\sum_{i=1}^{r_n} U_i / \sum_{j=1}^{n+1} U_j$ can be easily shown to follow $Beta(r_n, n+1-r_n)$. □

To finalize constructing the proposed tolerance region in (3.10), we need to identify a suitable $r_n$ which can satisfy Definitions 3.1 or 3.2. Following Theorem 3.1, this is equivalent to finding $r_n$ to meet the following criteria,

$$(3.11) \qquad P(Beta(r_n, n+1-r_n) \geq \beta) = \gamma$$

or

$$(3.12) \qquad E(Beta(r_n, n+1-r_n)) = \beta.$$

For (3.12), $r_n$ can be easily solved as

$$r_n = (n+1)\beta,$$

since $E(Beta(a,b)) = \frac{a}{a+b}$. For (3.11), it is not easy to find an analytical solution. Alternatively, we can obtain an approximation of the solution using the asymptotic result stated in Theorem 3.2 below.

REMARK 3.1 (Multivariate spacings as statistically equivalent blocks). For a multivariate sample of size $n$, Tukey (in [24]) considered a partition of the sample space into $n+1$ disjoint blocks as *statistically equivalent blocks* if the followings are satisfied:

(a) the coverages of the $(n+1)$ blocks add up to 1;
(b) the joint distribution of the coverages of the $(n+1)$ blocks are completely symmetrical;
(c) if the coverages of the $(n+1)$ blocks are taken as barycentric coordinates on an $n$-simplex, the distribution over the simplex is uniform;
(d) the sum of the coverages of any $k$ preselected blocks of the $(n+1)$ follows $Beta(k, n-k+1)$.

From the proof of Theorem 3.1, we can see that our multivariate spacings $\{MS_1, \ldots, MS_{n+1}\}$ satisfy the above conditions and can be viewed as statistically equivalent blocks. Note that the blocks as in our multivariate spacings are automatically determined by the given data. In contrast, the statistically equivalent blocks considered in [24] and other follow-ups all need to decide on the shape of the blocks before forming the blocks. Therefore, we view our multivariate spacings as an ideal data driven realization of statistically equivalent blocks. Moreover, the inherited property of data depth allows our statistical equivalent blocks to follow more closely the data structure and also be completely nonparametric.

THEOREM 3.2. *As $n \to \infty$, if $r_n$ satisfies*
$$\sqrt{n}\left(\frac{r_n}{n} - \beta\right) \to \xi_\gamma \sqrt{\beta(1-\beta)},$$
*where $\xi_\gamma$ is the $\gamma$-quantile of the standard normal distribution [i.e., $\Phi(\xi_\gamma) = \gamma$], then*
$$P(P_F(O_{Z^{[r_n]}}) \geq \beta) \to \gamma.$$

PROOF. Recall that $Y_i = P_F(X : D_F(X) > Z_i)$, $i = 1, \ldots, n$, with the order statistics $Y_{[1]} \leq \cdots \leq Y_{[n]}$. Let $\omega_n = \#\{i : Y_i < \beta\}$. Then we obtain $P(P_F(O_{Z^{[r_n]}}) \geq \beta) = P(\omega_n \leq r_n)$. Furthermore, since $Y_i$'s can be viewed as an i.i.d. sample from $U[0,1]$, $\omega_n$ follows the binomial distribution with parameter $(n, \beta)$. Therefore, as $n \to \infty$, we have
$$P(P_F(O_{Z^{[r_n]}}) \geq \beta) = P(\omega_n \leq r_n) \to \Phi\left(\frac{r_n - n\beta}{\sqrt{n\beta(1-\beta)}}\right) = \gamma. \qquad \square$$

3.2. *Properties of tolerance regions*: $F$ *is unknown.* If $F$ is unknown, the tolerance region is then constructed from the sample spacings in a similar fashion as in (3.10). More specifically, recall that $\hat{Z}_i = D_{F_n}(X_i)$, $i = 1, \ldots, n$, and that $\hat{Z}^{[1]} \geq \cdots \geq \hat{Z}^{[n]}$ are the descending estimated depth values corresponding to the depth order statistics $X_{[1]}, \ldots, X_{[n]}$. The tolerance region is then the union of a suitable number of the inner sample spacings. More precisely, the proposed tolerance region can be expressed as

$$(3.13) \qquad O^n_{\hat{Z}^{[r_n]}} = \bigcup_{i=1}^{r_n} \widehat{MS}_i \equiv \{X : D_{F_n}(X) > \hat{Z}^{[r_n]}\},$$

where $\widehat{MS}_i$ is the $i$th sample spacing, as defined in (2.4).

To establish the asymptotic properties for $O^n_{\hat{Z}^{[r_n]}}$ which are analogous to those for $O_{Z^{[r_n]}}$, we require the followings on the data depth $D_{F_n}(\cdot)$ used in the derivation of the spacings.

(i) If $F$ is absolutely continuous, then $D_{F_n}(x)$ is uniformly consistent almost surely, that is, as $n \to \infty$,

$$(3.14) \qquad d_n = \sup_x |D_{F_n}(x) - D_F(x)| \to 0 \quad \text{a.s.}$$

(ii) If $F$ is an elliptic distribution with the location-scatter parameter $(\mu, \Sigma)$ (i.e., its density assumes the form $f(x) = |\Sigma|^{-1/2} g((x-\mu)'\Sigma^{-1}(x-\mu))$), then its elliptic contour can be expressed as $e(x) = (x-\mu)'\Sigma^{-1}(x-\mu)$. In this case, the level sets (or contours) of $D_F(x)$ are also in the form of $\{x : e(x) = c\}$ for some $e(x)$. Furthermore, $D_F(x)$ is a strictly monotone function of $c$, which implies that for any $c > 0$,

$$(3.15) \qquad P(X : D_F(X) = c) = 0.$$

The discussion of (i) and (ii) under the simplicial depth can be found in [16]. Further discussions of depth contours can be found in [13, 17] and [31]. Under the assumptions (i) and (ii), we now present the main results of the section.

THEOREM 3.3. *Assume that conditions* (i) *and* (ii) *hold for the depth* $D_{F_n}(\cdot)$ *used in deriving the spacings. For any $\varepsilon > 0$, if the sequences* $r_{[n;1]}$, $r_{[n;2]}$ *and* $r_{[n;3]}$ $(1 \leq r_{[n;j]} \leq n, j = 1, 2, 3)$ *satisfy, as* $n \to \infty$,

$$\sqrt{n}\left(\frac{r_{[n;1]}}{n} - (\beta + \varepsilon)\right) \to \xi_\gamma \sqrt{\beta(1-\beta)},$$

$$\sqrt{n}\left(\frac{r_{[n;2]}}{n} - (\beta - \varepsilon)\right) \to \xi_\gamma \sqrt{\beta(1-\beta)},$$

$$\frac{r_{[n;3]}}{n+1} \to \beta,$$

*then*

$$\lim_{n \to \infty} P(P_F(O^n_{\hat{Z}^{[r_{[n;1]}]}}) \geq \beta) \geq \gamma,$$

$$\lim_{n \to \infty} P(P_F(O^n_{\hat{Z}^{[r_{[n;2]}]}}) \geq \beta) \leq \gamma$$

*and*

$$\lim_{n \to \infty} E(P_F(O^n_{\hat{Z}^{[r_{[n;3]}]}})) = \beta.$$

The proof of Theorem 3.3 is somewhat involved and is given in the Appendix.

REMARK 3.2. Since $\varepsilon$ in Theorem 3.3 can be arbitrarily small, to obtain $r_n$ satisfying $P(P_F(O^n_{\hat{Z}^{[r_n]}}) \geq \beta) = \gamma$, we may in practice simply take $\varepsilon = 0$ and calculate $r_n$ by solving

$$r_n = n\beta + \xi_\gamma \sqrt{n\beta(1-\beta)}.$$

If $r_n$ is not an integer, we use $\lfloor r_n \rfloor$ or $\lceil r_n \rceil$, depending on which of the following is closer to $\gamma$,

$$P(Beta(\lfloor r_n \rfloor, n+1-\lfloor r_n \rfloor) \geq \beta) \quad \text{and} \quad P(Beta(\lceil r_n \rceil, n+1-\lceil r_n \rceil) \geq \beta).$$

3.3. *Asymptotic minimum property.* So far, we have justified the proposed tolerance regions according to Definitions 3.1 and 3.2. Next, we will show that under a certain class of distributions (including elliptical distributions), the proposed tolerance regions are *asymptotically minimal*. This property is clearly desirable.

Asymptotically minimal tolerance regions were first considered in [3] by Chatterjee and Patra. Assume that the sample $X_1, X_2, \ldots, X_n$ is drawn from $F(\cdot) \in \Re^p$ which has a density function $f(\cdot)$. Let $\lambda(\cdot)$ denote the $p$-dimensional Lebesgue measure. Consider the set:

$$G_f(v) = P_F(f(X) \leq v).$$

Assume that all levels set of $f$ have Lebesgue measures zero, namely $\lambda\{x : f(x) = v\} = 0$ for any $v$. Chatterjee and Patra considered the following $\beta$-*content* tolerance region formed by density level sets:

(3.16) $$R_{f,\beta} = \{x : f(x) > \xi_{f,1-\beta}\},$$

where $\xi_{f,1-\beta}$ is the $(1-\beta)$-quantile of the random variable $f(X)$. In other words, $\xi_{f,1-\beta}$ is a solution of $G_f(v) = 1 - \beta$. It can be shown that

$$P_F(R_{f,\beta}) = \beta,$$

and that, among all subsets whose probability content with respect to $F$ is at least $\beta$, the subset $R_{f,\beta}$ is *minimal* in the sense of having the smallest Lebesgue measure.

DEFINITION 3.3. *A sequence of $\beta$-content tolerance regions $\{S_n\}$ is called asymptotically minimal if*

$$\lambda(S_n \Delta R_{f,\beta}) \xrightarrow{p} 0 \quad \text{as } n \to \infty.$$

Here $(A \Delta B)$ indicates the symmetric difference between sets $A$ and $B$.

In the finite sample case, [3] replaced $f$ with a density estimate to obtain a sample tolerance region, and showed the asymptotic minimum property. Clearly, the quality of the obtained tolerance region depends on the density estimation approach used.

Under the approach with data depth $D(\cdot)$, we consider

$$G_D(v) = P_F(D_F(X) \leq v).$$

Denote by $\eta_{1-\beta}$ the $(1-\beta)$-quantile of the random variable $D_F(X)$. In other words,

$$G_D(\eta_{1-\beta}) = P_F(D_F(X) \leq \eta_{1-\beta}) = 1 - \beta.$$

Clearly,

(3.17) $$R_{D,\beta} = \{X : D_F(X) > \eta_{1-\beta}\},$$

is the true depth-based $\beta$-content tolerance region. Definition 3.3 can then be modified for the approach using depth $D(\cdot)$ as:

DEFINITION 3.4. *A sequence of $\beta$-content tolerance regions $\{S_n\}$ is called asymptotically minimal w.r.t. the depth function $D(\cdot)$ if*

$$\lambda(S_n \Delta R_{D,\beta}) \xrightarrow{p} 0 \quad \text{as } n \to \infty.$$

In the next two theorems we show that our proposed tolerance regions $O_{Z^{[r_n]}}$ and $O^n_{\hat{Z}^{[r_n]}}$ are asymptotically minimal w.r.t. the chosen depth.

THEOREM 3.4. *If condition (3.15) holds for the underlying depth $D(\cdot)$, $O_{Z^{[r_n]}}$ and $O^n_{\hat{Z}^{[r_n]}}$ are asymptotically minimal w.r.t. $D(\cdot)$. Specifically, for $r_n$ satisfying*

$$\sqrt{n}\left(\frac{r_n}{n} - \beta\right) \to \xi_\gamma \sqrt{\beta(1-\beta)},$$

*we have*

$$\lambda(O_{Z^{[r_n]}} \Delta R_{D,\beta}) \xrightarrow{p} 0 \quad \text{and} \quad \lambda(O^n_{\hat{Z}^{[r_n]}} \Delta R_{D,\beta}) \xrightarrow{p} 0.$$

The proof of this theorem is given in the Appendix.

Note that, for elliptical distributions, condition (3.15) holds for all the depth notions mentioned in Section 2.1, and thus

$$R_{D,\beta} = \{X : D_F(X) > \eta_{1-\beta}\} = \{x : f(x) > \xi_{f,1-\beta}\} = R_{f,\beta}.$$

Consequently, Theorem 3.4 leads to the corollary below which implies that our proposed tolerance regions are asymptotically minimal under elliptical distributions.

COROLLARY 3.1. *For elliptical distributions, we have, under condition (3.15), as $n \to 0$,*

$$\lambda(O_{Z^{[r_n]}} \Delta R_{f,\beta}) \xrightarrow{p} 0,$$

$$\lambda(O^n_{\hat{Z}^{[r_n]}} \Delta R_{f,\beta}) \xrightarrow{p} 0.$$

**4. Simulation and comparison studies.** In this section, we present some simulation studies to illustrate the performance of our tolerance regions $O_{Z^{[r_n]}}$ and $O^n_{\hat{Z}^{[r_n]}}$. The simulation procedure is outlined in the following steps. Assume that $F$ is absolutely continuous.

Step 1. Generate a random sample $\{X_1, X_2, \ldots, X_n\}$ from $F$. Calculate the depths of $X_i$'s with respect to the data cloud and identify the $r_n$th largest depth (i.e. $Z^{[r_n]}$) as the threshold for forming our tolerance region, where

$$r_n = \begin{cases} n\beta + \sqrt{n}\xi_\gamma\sqrt{\beta(1-\beta)}, & \text{if (3.1) is required,} \\ (n+1)\beta, & \text{if (3.2) is required.} \end{cases}$$

Step 2. Generate another random sample, $\{X_1^*, X_2^*, \ldots, X_n^*\}$, from $F$. Calculate the depth of $X_i^*$ with respect to the original samples, $\{X_1, X_2, \ldots, X_n\}$, and obtain the proportion of $X_i^*$'s which assume depth value greater than the threshold obtained in Step 1. This proportion is denoted as $\tilde{\beta}$.

Step 3. Repeat Step 2 $m$ times and use the average of the $m$ $\tilde{\beta}$'s obtained in this manner as an approximate of $P_F(O^n_{\hat{Z}^{[r_n]}})$. Denote this average by $\bar{\beta}$. If a $\beta$-*content* tolerance region in (3.1) is sought, we check whether or not $\bar{\beta} > \beta$. Let $I_{\{\bar{\beta}>\beta\}}$ be 1 if the event $\{\bar{\beta} > \beta\}$ occurs, and 0 otherwise. If a $\beta$-*expectation* tolerance region in (3.2) is sought, we simply record the $\bar{\beta}$.

Step 4. Repeat Steps 1–3 sufficiently many, say $M$, times. For a $\beta$-*content* tolerance region, we estimate the confidence level $\gamma$ by $\hat{\gamma} = \sum_{i=1}^M I_{\{\bar{\beta}>\beta\}}/M$. For a $\beta$-*expectation* tolerance region, we estimate $\beta$ simply by the average of the $M$ $\bar{\beta}$'s, namely $\hat{\beta} = \sum_{i=1}^M \bar{\beta}_i/M$.

Throughout our simulation study, we set $\beta = 90\%$, $\gamma = 95\%$, $m = 100$, $M = 1000$ and $n = 300$ and 1000. The simplicial depth is used to calculate all depth values unless specified otherwise. The results are presented in Table 1. From Table 1, we can see that all the estimates are fairly close to the nominal levels. We also present in Table 1, within brackets, the simulation results using the approach given in [6] by Di Bucchianico, Einmahl and Mushkudiani. The coverage from this approach is consistently lower than the nominal value, and also, generally speaking, the difference between the achieved and nominal coverage is larger than that of ours. Thus, our proposed tolerance region is better in terms of achieving the desired tolerance level. Moreover, as observed in [6], if the dimension of the underlying distribution increases, the approach there needs an adjustment to reflect such a change in the target nominal value to prevent the achieved coverage from falling too much below. The adjustment suggested in [6] seems somewhat ad hoc. Adding all these observations to the fact that our approach does not require additional assumptions (e.g., shape of the tolerance region), our approach clearly yields more favorable nonparametric multivariate tolerance regions.

To help visualize the outcome of our constructions, we present further in Figure 3 our proposed tolerance region for bivariate normal and exponential distributions. Under each distribution, the sample size is 500 and the tolerance regions shown are aiming for nominal values $\beta = 90\%$ and $\gamma = 95\%$. Note that since there is no explicit formula for the simplicial depth, there is no explicit expression for the proposed tolerance region $O_{\hat{Z}^{[r_n]}}^n = \{X : D_{F_n}(X) > \hat{Z}^{[r_n]}\}$. (Here $r_n$ is to be determined according to Remark 3.2.) In practice, we can simply present the convex hull spanning all the sample points which achieve higher depth value than $\hat{Z}^{[r_n]}$ as the estimated tolerance region. The algorithm provided in [1] can help determine such convex hulls.

As seen from the plots, our tolerance regions can capture the underlying geometric shapes of the data. For the bivariate normal distribution, the

TABLE 1
*The achieved confidence levels of the proposed $\beta$-content and $\beta$-expectation tolerance regions when $\gamma = 95\%$ and $\beta = 90\%$*

|  | $\hat{\gamma}$ | | $\hat{\beta}$ | |
|---|---|---|---|---|
| $F$ | $n = 300$ | $n = 1000$ | $n = 300$ | $n = 1000$ |
| Bivariate Normal | 0.954 | 0.949 | 0.90131 [0.877] | 0.90005 [0.887] |
| Bivariate Cauchy | 0.963 | 0.961 | 0.90036 [0.862] | 0.90061 [0.863] |
| Bivariate Exponential | 0.941 | 0.943 | 0.90043 [0.885] | 0.89985 [0.890] |

Results in [] are obtained using approach in [6].

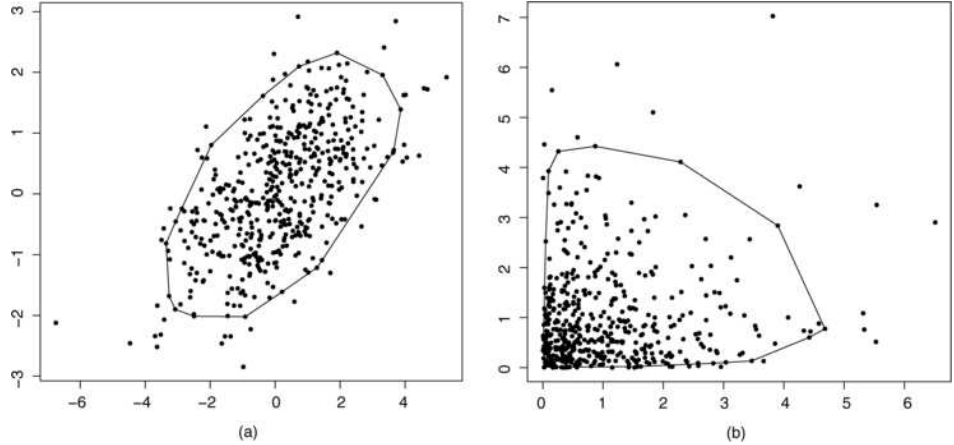

FIG. 3.  *Tolerance region for:* (a) *a bivariate normal sample;* (b) *a bivariate exponential sample.*

region has the elliptical shape. For the bivariate exponential distribution, the region has a triangular shape fanning upper-right. Overall, our tolerance region focuses more on the central part of the data and also follows the expansion of the probability mass. For example, for the bivariate exponential distribution, our tolerance region does not include the observations near the origin, since their positions are relatively outlying with respect to the center of the distribution. In contrast, the tolerance region obtained by using the method proposed in [6] must include those observations since they have high density. However, these observations may not be acceptable in practice, since they may be considered extreme according to the underlying distribution. This point can be made more pronounced by incurring a small perturbation to create a thinning gap between the points near the origin and the rest of the data. Therefore, the design of our tolerance region in representing a central region of the data is naturally built to prevent the region from including observations which are likely to be extreme.

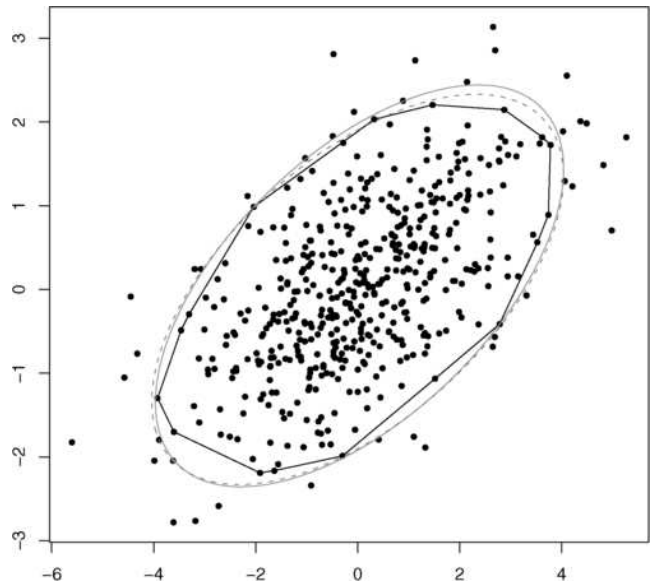

FIG. 4.  *Tolerance regions for a bivariate normal sample: unknown $F$ vs. known $F$.*

As discussed in Section 3.1, when the underlying distribution $F$ is known, we can construct the tolerance region which satisfies exactly the preset requirement in (3.1) or (3.2). When $F$ is unknown, we propose a method for constructing the tolerance region based on the sample only and develop their asymptotic properties. From the asymptotic point of view, the proposed method also satisfies the preset requirement. To assess the performance of this proposed tolerance region in the setting of a finite sample size, we conduct another simulation study. In the same bivariate normal setting as above, we use the Mahalanobis depth to construct the tolerance regions separately under $F$ is known and $F$ is unknown. One advantage in using the Mahalanobis depth is that it has a closed form for both population and sample versions, and thus we can obtain the exact proposed tolerance regions. Figure 4 shows the two tolerance regions. The dashed ellipse is the one when $F$ is known, which is the true tolerance region. The solid ellipse is the proposed tolerance region when $F$ is unknown. The two regions almost coincide with each other, which clearly implies that the finite sample performance of the proposed tolerance region is quite satisfactory. The convex hull in Figure 4 is formed by only the sample points which have higher depth value than $\hat{Z}^{[r_n]}$. (This generally reduces tremendously the computational effort, and is more practical, especially when the depth itself does not have a simple closed form.) It is not surprising that the convex hull is located inside the solid ellipse. Note that difference between these two regions is not significant. Therefore, although the convex hull formed by those central points is not the exact proposed tolerance region, it is presented here to show that it can be a viable alternative that provides a practical solution. To use the tolerance region in terms of certifying specifications, we determine whether a new observation is in the tolerance region by first calculating its depth with respect to the given sample and then simply comparing the depth value to the threshold $\hat{\tilde{Z}}^{[r_n]}$. This is a relatively straightforward task in practice.

**5. Concluding remarks and future research.** In this paper, we introduced multivariate spacings based on the ordering derived from data depth.

They satisfy all properties one would expect of a notion of spacings. Moreover, they are nonparametric and they reflect well the geometry of the underlying distribution. We show how to use these spacings to construct tolerance regions for multivariate distributions. The construction of our tolerance regions can be viewed as a multi-dimensional generalization of the Wilks' method.

Given that our spacings are derived from data depth, the resulting tolerance regions are always connected and naturally located in the "central" region of the data set. This is an important property in applications: in practice, specifications of products are not given in disconnected patches and a single production line is generally designed to produce continuous measurement around a target value. The connected tolerance region ensures that the capability of production processes can be achieved. This point was also discussed and illustrated with Figure 3 in Section 4.

One important direction for the applications of univariate spacings is nonparametric inference. It includes many existing rank tests and goodness-of-fit tests. A survey of these tests as well as relevant references can be found in [22]. In forthcoming papers, we shall explore our multivariate spacings in the development of nonparametric inference. We generalize many of the existing approaches on univariate spacings. We also study multivariate distributional characterizations using our multivariate spacings.

## APPENDIX

PROOF OF THEOREM 3.3. To prove Theorem 3.3, we need the following two lemmas.

LEMMA A.1. *For any $r$, $|\hat{Z}^{[r]} - Z^{[r]}| \leq d_n = \sup_x |D_{F_n}(x) - D_F(x)| \to 0$, a.s.*

PROOF. From the definition of $d_n$, we have
$$|\hat{Z}_i - Z_i| = |D_{F_n}(X_i) - D_F(X_i)| \leq d_n, \qquad i = 1, \ldots, n.$$
Then for any $r$,
$$\#\{i : Z_i < \hat{Z}^{[r]} - d_n\} \leq \#\{i : \hat{Z}_i < \hat{Z}^{[r]}\} < r.$$
Therefore, $\hat{Z}^{[r]} - d_n \leq Z^{[r]}$. Similarly, we can show that $Z^{[r]} - d_n \leq \hat{Z}^{[r]}$. The claim of the lemma thus follows. □

LEMMA A.2. *Suppose that $a_n$ and $b_n$ are two sequences of random variables such that for some random variable $a$ taking values on $[0, \infty]$, and* $a_n \to a$ and $b_n \to a$ on a positive measure subset of the sample space, say $S$, as $n \to \infty$. Then under the assumptions (i)–(ii),

$$P_F\{O^n_{a_n} \Delta O_{b_n}\} \to 0 \qquad \text{a.s. on the set } S,$$

where $O^n_{a_n} = \{x : D_{F_n}(x) \geq a_n\}$, $O_{b_n} = \{x : D_F(x) \geq b_n\}$, and $A\Delta B = (A \cup B) \setminus (A \cap B)$.

Lemma A.2 with proof is given in [13].

We now proceed with the proof of Theorem 3.3. Assume that $\frac{r_n}{n} \to q$, as $n \to \infty$, then the consistency property of a sample quantile shows that $Z^{[r_n]} \to \eta_q$, a.s., where $\eta_q$ is the upper $q$th quantile of $Z = D_F(X)$. Clearly, Lemma A.1 immediately implies $\hat{Z}^{[r_n]} \to \eta_q$, a.s. Let $a_n = \hat{Z}^{[r_n]}$ and $b_n = Z^{[r_n]}$. Following Lemma A.2, we have

(A.1) $$P_F(O^n_{\hat{Z}^{[r_n]}} \Delta O_{Z^{[r_n]}}) \to 0 \qquad \text{a.s.}$$

Denote $A_n = P_F(O^n_{\hat{Z}^{[r_n]}})$, $B_n = P_F(O^n_{\hat{Z}^{[r_n]}} \Delta O_{Z^{[r_n]}})$ and $C_n = P_F(O_{Z^{[r_n]}})$. Then

$$\begin{aligned} P(C_n \geq \beta) &\leq P(A_n + B_n \geq \beta) \\ &= P(A_n + B_n \geq \beta \cap B_n \leq \varepsilon) + P(A_n + B_n \geq \beta \cap B_n > \varepsilon) \qquad \forall \varepsilon > 0 \\ &\leq P(A_n \geq \beta - \varepsilon) + P(B_n > \varepsilon). \end{aligned}$$

From (A.1), we have, as $n \to \infty$,

$$P(B_n > \varepsilon) \to 0$$

and

$$\lim_{n \to \infty} P(C_n \geq \beta) \leq \lim_{n \to \infty} P(A_n \geq \beta - \varepsilon) \qquad \forall \varepsilon > 0.$$

Therefore,

$$\lim_{n \to \infty} P(A_n \geq \beta) \geq \lim_{n \to \infty} P(C_n \geq \beta + \varepsilon) \qquad \forall \varepsilon > 0.$$

If $r_n$ satisfies

$$\sqrt{n}\left(\frac{r_n}{n} - (\beta + \varepsilon)\right) \to \xi_\gamma \sqrt{\beta(1-\beta)} \quad \text{and} \quad \lim_{n \to \infty} P(C_n \geq \beta + \varepsilon) \to \gamma,$$

then

$$\lim_{n \to \infty} P(A_n \geq \beta) \geq \gamma.$$

Similarly, we have

$$P(A_n \geq \beta) \leq P(C_n + B_n \geq \beta) \leq P(C_n \geq \beta - \varepsilon) + P(B_n > \varepsilon),$$

and thus
$$\lim_{n\to\infty} P(A_n \geq \beta) \leq \lim_{n\to\infty} P(C_n \geq \beta - \varepsilon) \qquad \forall \varepsilon > 0.$$

If $r_n$ satisfies
$$\sqrt{n}\left(\frac{r_n}{n} - (\beta - \varepsilon)\right) \to \xi_\gamma \sqrt{\beta(1-\beta)} \quad \text{and} \quad \lim_{n\to\infty} P(C_n \geq \beta - \varepsilon) \to \gamma,$$

then $\lim_{n\to\infty} P(A_n \geq \beta) \leq \gamma$.

Regarding the $\beta$-expectation tolerance region, we now have, following (A.1),
$$B_n \to 0 \qquad \text{a.s. as } n \to \infty.$$

Since $B_n \leq 1$, $B_n$ is uniformly integrable. Thus $\lim_{n\to\infty} E(B_n) = 0$. Since $E(A_n) \leq E(C_n + B_n) \leq E(C_n) + E(B_n)$, if $r_n$ satisfies that $\frac{r_n}{n+1} \to \beta$, then $\lim_{n\to\infty} E(C_n) = \beta$, and

$$\limsup_{n\to\infty} E(A_n) \leq \beta. \tag{A.2}$$

Similarly, $E(C_n) \leq E(A_n + B_n) \leq E(A_n) + E(B_n)$. Then we have

$$\beta \leq \liminf_{n\to\infty} E(A_n). \tag{A.3}$$

Combining (A.2) and (A.3), we obtain $\lim_{n\to\infty} E(A_n) = \beta$, and hence the proof of Theorem 3.3. □

PROOF OF THEOREM 3.4. Since $\frac{r_n}{n} \to \beta$,
$$G_D(Z^{[r_n]}) = P_F(D_F(X) \leq Z^{[r_n]}) \to 1 - \beta$$

which implies

$$Z^{[r_n]} \xrightarrow{p} \eta_{1-\beta}. \tag{A.4}$$

Moreover, the following
$$O_{Z^{[r_n]}} \Delta R_{D,\beta}$$
$$= \{X : D_F(X) > Z^{[r_n]}, D_F(X) \leq \eta_{1-\beta}\}$$
$$\cup \{X : D_F(X) \leq Z^{[r_n]}, D_F(X) > \eta_{1-\beta}\}$$
$$\subset \{X : \eta_{1-\beta} - |Z^{[r_n]} - \eta_{1-\beta}| \leq D_F(X) \leq \eta_{1-\beta} + |Z^{[r_n]} - \eta_{1-\beta}|\}$$

immediately implies that

$$\lambda(O_{Z^{[r_n]}} \Delta R_{D,\beta})$$
$$\leq \lambda\{X : \eta_{1-\beta} - |Z^{[r_n]} - \eta_{1-\beta}| \leq D_F(X) \leq \eta_{1-\beta} + |Z^{[r_n]} - \eta_{1-\beta}|\} \tag{A.5}$$
$$= \delta(|Z^{[r_n]} - \eta_{1-\beta}|),$$

where $\delta(u) = \lambda\{X : \eta_{1-\beta} - u \leq D_F(X) \leq \eta_{1-\beta} + u\}$. By assumption (3.15), $\delta(0) = \lambda\{X : D_F(X) = \eta_{1-\beta}\} = 0$. Also $\delta(u)$ is right continuous at 0 because of the continuity of $D_F(x)$. Therefore, (A.4) implies that $\delta(|Z^{[r_n]} - \eta_{1-\beta}|) \xrightarrow{p} 0$. Following (A.5), we finally obtain $\lambda(O_{Z^{[r_n]}} \Delta R_{D,\beta}) \xrightarrow{p} 0$.

The rest of the proof regarding $O^n_{\hat{Z}^{[r_n]}}$ can be derived similarly. Following Lemma A.1 and the definition of $d_n$, we obtain

$$\begin{aligned}
O^n_{\hat{Z}^{[r_n]}} &\Delta R_{D,\beta} \\
&= \{X : D_{F_n}(X) > \hat{Z}^{[r_n]}, D_F(X) \leq \eta_{1-\beta}\} \\
&\quad \cup \{X : D_{F_n}(X) \leq \hat{Z}^{[r_n]}, D_F(X) > \eta_{1-\beta}\} \\
&\subset \{X : D_F(X) > Z^{[r_n]} - 2d_n, D_F(X) \leq \eta_{1-\beta}\} \\
&\quad \cup \{X : D_F(X) \leq Z^{[r_n]} + 2d_n, D_F(X) > \eta_{1-\beta}\} \\
&\subset \{X : \eta_{1-\beta} - |Z^{[r_n]} - \eta_{1-\beta}| - 2d_n \\
&\qquad \leq D_F(X) \leq \eta_{1-\beta} + |Z^{[r_n]} - \eta_{1-\beta}| + 2d_n\}.
\end{aligned}$$

Finally, since $Z^{[r_n]} \xrightarrow{p} \eta_{1-\beta}$ and $d_n \xrightarrow{p} 0$, we have

$$\begin{aligned}
\lambda(O^n_{\hat{Z}^{[r_n]}} &\Delta R_{D,\beta}) \\
&\leq \lambda\{X : \eta_{1-\beta} - |Z^{[r_n]} - \eta_{1-\beta}| - 2d_n \\
&\qquad \leq D_F(X) \leq \eta_{1-\beta} + |Z^{[r_n]} - \eta_{1-\beta}| + 2d_n\} \\
&= \delta(|Z^{[r_n]} - \eta_{1-\beta}| + 2d_n) \\
&\xrightarrow{p} 0.
\end{aligned}$$

This completes the proof. □

**Acknowledgment.** Jun Li would also like to acknowledge the graduate assistantship provided by the Department of Statistics, Rutgers University.

DEPARTMENT OF STATISTICS  
UNIVERSITY OF CALIFORNIA  
RIVERSIDE, CALIFORNIA 92521-0138  
USA  
E-MAIL: jun.li@ucr.edu

DEPARTMENT OF STATISTICS  
RUTGERS UNIVERSITY, HILL CENTER  
PISACATAWAY, NEW JERSEY 08854-8019  
USA  
E-MAIL: rliu@stat.rutgers.edu